\documentclass[12pt,reqno]{amsart}

\usepackage[utf8]{inputenc}
\usepackage{latexsym} 
\usepackage{amsfonts}
\usepackage{amssymb} 
\usepackage{amsmath}
\usepackage{enumerate}
\usepackage{hyperref}  

\usepackage{bm}
\usepackage{amsfonts}

\usepackage{graphicx}
\usepackage{eufrak}
\usepackage{pstricks}
\usepackage{graphicx}
\usepackage{dsfont}

\begin{document}

\newtheorem{theorem}{Theorem}[section]
\newtheorem{problem}{Problem} [section]
\newtheorem{definition}{Definition} [section]
\newtheorem{lemma}{Lemma}[section]
\newtheorem{proposition}{Proposition}[section]
\newtheorem{corollary}{Corollary}[section]
\newtheorem{example}{Example}[section]
\newtheorem{conjecture}{Conjecture} 
\newtheorem{algorithm}{Algorithm} 
\newtheorem{exercise}{Exercise}[section]
\newtheorem{remarkk}{Remark}[subsection]
 
\newcommand{\be}{\begin{equation}} 
\newcommand{\ee}{\end{equation}} 
\newcommand{\bea}{\begin{eqnarray}} 
\newcommand{\eea}{\end{eqnarray}} 

\newcommand{\eeq}{\end{equation}} 

\newcommand{\eeqn}{\end{eqnarray}} 
\newcommand{\beaa}{\begin{eqnarray*}} 
\newcommand{\eeaa}{\end{eqnarray*}}

\newcommand{\lip}{\langle} 
\newcommand{\rip}{\rangle}

\newcommand{\uu}{\underline} 
\newcommand{\oo}{\overline} 
\newcommand{\La}{\Lambda} 
\newcommand{\la}{\lambda} 
\newcommand{\eps}{\varepsilon} 
\newcommand{\om}{\omega} 
\newcommand{\Om}{\Omega} 
\newcommand{\ga}{\gamma} 
\newcommand{\rrr}{{\Bigr )}} 
\newcommand{\qqq}{{\Bigl\|}} 
 
\newcommand{\dint}{\displaystyle\int} 
\newcommand{\dsum}{\displaystyle\sum} 
\newcommand{\dfr}{\displaystyle\frac} 
\newcommand{\bige}{\mbox{\Large\it e}} 
\newcommand{\integers}{{\Bbb Z}} 
\newcommand{\rationals}{{\Bbb Q}} 
\newcommand{\reals}{{\rm I\!R}} 
\newcommand{\realsd}{\reals^d} 
\newcommand{\realsn}{\reals^n} 
\newcommand{\NN}{{\rm I\!N}} 
\newcommand{\DD}{{\rm I\!D}} 
\newcommand{\degree}{{\scriptscriptstyle \circ }} 
\newcommand{\dfn}{\stackrel{\triangle}{=}} 
\def\complex{\mathop{\raise .45ex\hbox{${\bf\scriptstyle{|}}$} 
     \kern -0.40em {\rm \textstyle{C}}}\nolimits} 
\def\hilbert{\mathop{\raise .21ex\hbox{$\bigcirc$}}\kern -1.005em {\rm\textstyle{H}}} %Hilbert space 
\newcommand{\RAISE}{{\:\raisebox{.6ex}{$\scriptstyle{>}$}\raisebox{-.3ex} 
           {$\scriptstyle{\!\!\!\!\!<}\:$}}} % >< one above each other 
 
\newcommand{\hh}{{\:\raisebox{1.8ex}{$\scriptstyle{\degree}$}\raisebox{.0ex} 
           {$\textstyle{\!\!\!\! H}$}}} 

\newcommand{\OO}{\won} 
\newcommand{\calA}{{\mathcal A}} 
\newcommand{\calB}{{\cal B}} 
\newcommand{\calC}{{\cal C}} 
\newcommand{\calD}{{\cal D}} 
\newcommand{\calE}{{\cal E}} 
\newcommand{\calF}{{\mathcal F}} 
\newcommand{\calG}{{\cal G}} 
\newcommand{\calH}{{\cal H}} 
\newcommand{\calK}{{\cal K}} 
\newcommand{\calL}{{\mathcal L}} 
\newcommand{\calM}{{\mathcal M}} 
\newcommand{\calO}{{\cal O}} 
\newcommand{\calP}{{\cal P}} 
\newcommand{\calU}{{\mathcal U}} 
\newcommand{\calX}{{\cal X}} 
\newcommand{\calXX}{{\cal X\mbox{\raisebox{.3ex}{$\!\!\!\!\!-$}}}} 
\newcommand{\calXXX}{{\cal X\!\!\!\!\!-}} 
\newcommand{\gi}{{\raisebox{.0ex}{$\scriptscriptstyle{\cal X}$} 
\raisebox{.1ex} {$\scriptstyle{\!\!\!\!-}\:$}}} 
\newcommand{\intsim}{\int_0^1\!\!\!\!\!\!\!\!\!\sim} 
\newcommand{\intsimt}{\int_0^t\!\!\!\!\!\!\!\!\!\sim} 
\newcommand{\pp}{{\partial}} 
\newcommand{\al}{{\alpha}} 
\newcommand{\sB}{{\cal B}} 
\newcommand{\sL}{{\cal L}} 
\newcommand{\sF}{{\cal F}} 
\newcommand{\sE}{{\cal E}} 
\newcommand{\sX}{{\cal X}} 
\newcommand{\R}{{\rm I\!R}} 
\renewcommand{\L}{{\rm I\!L}} 
\newcommand{\vp}{\varphi} 
\newcommand{\N}{{\rm I\!N}} 
\def\ooo{\lip} 
\def\ccc{\rip} 
\newcommand{\ot}{\hat\otimes} 
\newcommand{\rP}{{\Bbb P}} 
\newcommand{\bfcdot}{{\mbox{\boldmath$\cdot$}}} 
 
\renewcommand{\varrho}{{\ell}} 
\newcommand{\dett}{{\textstyle{\det_2}}} 
\newcommand{\sign}{{\mbox{\rm sign}}} 
\newcommand{\TE}{{\rm TE}} 
\newcommand{\TA}{{\rm TA}} 
\newcommand{\E}{{\rm E\, }} 
\newcommand{\won}{{\mbox{\bf 1}}} 
\newcommand{\Lebn}{{\rm Leb}_n} 
\newcommand{\Prob}{{\rm Prob\, }} 
\newcommand{\sinc}{{\rm sinc\, }} 
\newcommand{\ctg}{{\rm ctg\, }} 
\newcommand{\loc}{{\rm loc}} 
\newcommand{\trace}{{\, \, \rm trace\, \, }} 
\newcommand{\Dom}{{\rm Dom}} 
\newcommand{\ifff}{\mbox{\ if and only if\ }} 
\newcommand{\nproof}{\noindent {\bf Proof:\ }} 
\newcommand{\nproofYWN}{\noindent {\bf Proof of Theorem~\cite{YWN}:\ }} 
\newcommand{\remark}{\noindent {\bf Remark:\ }} 
\newcommand{\remarks}{\noindent {\bf Remarks:\ }} 
\newcommand{\note}{\noindent {\bf Note:\ }} 
 \newcommand{\examples}{\noindent {\bf Examples:\ }} 
 
\newcommand{\boldx}{{\bf x}} 
\newcommand{\boldX}{{\bf X}} 
\newcommand{\boldy}{{\bf y}} 
\newcommand{\boldR}{{\bf R}} 
\newcommand{\uux}{\uu{x}} 
\newcommand{\uuY}{\uu{Y}} 
 
\newcommand{\limn}{\lim_{n \rightarrow \infty}} 
\newcommand{\limN}{\lim_{N \rightarrow \infty}} 
\newcommand{\limr}{\lim_{r \rightarrow \infty}} 
\newcommand{\limd}{\lim_{\delta \rightarrow \infty}} 
\newcommand{\limM}{\lim_{M \rightarrow \infty}} 
\newcommand{\limsupn}{\limsup_{n \rightarrow \infty}} 
 
\newcommand{\ra}{ \rightarrow } 

 \newcommand{\mlim}{\lim_{m \rightarrow \infty}}  
 \newcommand{\limm}{\lim_{m \rightarrow \infty}}  
 \newcommand{\nlim}{\lim_{n \rightarrow \infty}} 
%\newcommand{\imii}{\int_{-\infty}^{\infty}} 
%\newcommand{\imix}{\int_{-\infty}^x} 
%\newcommand{\ioi}{\int_o^\infty} 
 
%\newcommand{\ARROW}[1\right] 
  %{\begin{array}[t\right]{c}  \longrightarrow \\[-0.2cm\right] \textstyle{#1} \end{array} } 
 
%\newcommand{\AR} 
 %{\begin{array}[t\right]{c} 
  %\longrightarrow \\[-0.3cm\right] 
  %\scriptstyle {n\rightarrow \infty} 
 % \end{array}} 
 
%\newcommand{\pile}[2\right] 
  %{\left( \begin{array}{c}  {#1}\\[-0.2cm\right] {#2} \end{array} \right ) } 
 
%\newcommand{\floor}[1\right]{\left\lfloor #1 \right\rfloor} 
 
%for doing boldface subscripts etc., e.g. $G_{\mmbox{\boldx}}$ 
%\newcommand{\mmbox}[1\right]{\mbox{\scriptsize{#1}}} 
 
%fraction with round brackets 
%\newcommand{\ffrac}[2\right] 
 % {\left( \frac{#1}{#2} \right )} 

\newcommand{\one}{\frac{1}{n}\:} 
\newcommand{\half}{\frac{1}{2}\:} 
 
\def\le{\leq} 
\def\ge{\geq} 
\def\lt{<} 
\def\gt{>} 

%qed 
\def\squarebox#1{\hbox to #1{\hfill\vbox to #1{\vfill}}} 
\newcommand{\nqed}{\hspace*{\fill} 
           \vbox{\hrule\hbox{\vrule\squarebox{.667em}\vrule}\hrule}\bigskip}

\title[adapted coupling of real valued random variables]{A mass transport approach to the optimization of adapted couplings of real valued random variables.}

%\subtitle{Using  the  LaTex Template}

\author{R\'{e}mi Lassalle}

\address{Universit\'{e} Paris 9 (Dauphine), PSL, Place du Mar\'echal De Lattre De Tassigny, 75775 Paris Cedex 16, France}
\email{lassalle@ceremade.dauphine.fr}

\maketitle

\begin{abstract} 
In this work, we investigate an optimization problem over adapted couplings  between pairs of real valued random variables, possibly describing random times.  We relate those couplings to a specific class of  causal transport plans between probabilities on the set of real numbers endowed with a filtration, for which their provide a specific representation, which is motivated by a precise characterization of the corresponding deterministic transport plans. From this, under mild hypothesis, the existence of a solution   to the problem investigated here is obtained. Furthermore, several examples are provided, within this explicit framework.

\end{abstract}

%All acknowledgements should be placed in the back of the paper after Conclusions..

%\tableofcontents

\section{Introduction}

 For the sake of clarity, in this work, the set  $\mathbb{T}$ which within models encountered here may label times, will be the whole real line $\mathbb{R}$. We investigate the connections between adapted mass transport problems and optimization over fixed marginals couplings $(X,Y, (\Omega,\mathcal{A},\mathbb{P}))$ of the specific form :

\begin{definition} \label{DEFCA1}
Let $X :\Omega\to \mathbb{R}$ and $Y: \Omega \to \mathbb{R}$ be two random variables defined on a same $\mathbb{P}-$ complete probability space $(\Omega,\mathcal{A},\mathbb{P})$. We say that $(X,Y, (\Omega,\mathcal{A},\mathbb{P}))$ is an adapted coupling, if and only if, there exist two  random variables $\tau :\Omega \to \mathbb{R}$, and $Z :\Omega \to \mathbb{R}$, such that \begin{equation} Y= (X+Z) \mathds{1}_{\{X\leq \tau\}} + \tau \mathds{1}_{\{X>\tau\}}, \  \mathbb{P}-a.s., \label{Ydefcoupl} \end{equation} where $Z$ and $\tau$ match with the following conditions :
\begin{enumerate}
\item $Z\geq 0$, $\mathbb{P}-a.s.$
\item For any $t\in \mathbb{\mathbb{T}}$, such that $\mathbb{P}(\{X\geq t\}) >0$,  the condition $$\mathbb{P}(\{\tau\in A\}\cap \{X\in B\} | \{X\geq t\}) = \mathbb{P}(\{\tau\in A\} | \{X\geq t\})\mathbb{P}( \{X\in B\} | \{X\geq t\}),$$ holds, $\forall A,B\in\mathcal{B}(\mathbb{R})$ \textsc{Borel} sets of $\mathbb{R}$, such that $A\subset ]-\infty, t[$ and $B\subset [t,+\infty[$.
\end{enumerate}
Given $\eta,\nu \in M_1(\mathbb{R})$, two \textsc{Borel} probability measures on $\mathbb{R}$, we denote by $Cpl_a(\eta,\nu)$ the set of such adapted couplings  $(X,Y, (\Omega,\mathcal{A},\mathbb{P}))$ which further satisfy $p_X=\eta$ and $p_Y= \nu$, where $p_X$ (resp. $p_Y$) denotes the probability law of the random variable $X$ (resp. $Y$) on the probability space  $(\Omega,\mathcal{A},\mathbb{P})$. 
\end{definition}

In models, we may interpret $X$ as a time where a signal, for instance a postal letter, or a phone call, or any physical signal within experiments, is received, and $Y$ as the time where a response signal is emitted, while we interpret $\tau$ as a time waiting for the signal actually received at $X$, and $Z$ as a delay time of answer ; if the entry signal has not been received before the waiting time $\tau$ has been exhausted, it results into a signal sent at $\tau$ ; otherwise, a signal is emitted with a delay $Z$ after the entry signal has been received, before $\tau$,  at a time $X$. Thus, we interpret $(1)$ as a \textit{delay} hypothesis, while we interpret $(2)$ as a \textit{waiting time} hypothesis. Within this perspective, notice that a sufficient condition for $(2)$ is that  $\tau$ is $\mathbb{P}-$ independent to $X$, while another sufficient condition is obtained by taking $\tau=X$, $\mathbb{P}-a.s.$, whenever $Y\geq X$ holds $\mathbb{P}-$ almost surely. Therefore, taking now for granted the existence of such couplings, let $X$, $\tau$, $Z$  be three  random variables defined on a same $\mathbb{P}-$ complete probability space $(\Omega, \mathcal{A},\mathbb{P})$, where $Z$ is non-negative, where $\tau$ satisfies $(2)$, and define $Y$ by~(\ref{Ydefcoupl}) ; then $(X,Y, (\Omega,\mathcal{A},\mathbb{P}))$ is an adapted coupling according to Definition~\ref{DEFCA1}.  In particular, if $Z=0$, then in the previous example, we obtain $Y = \min(X, \tau) : $ thus $\tau$ suggests a particular case of so-called \textit{censoring times} (see \cite{PROTTER} p.124), which might involve further interpretations for $\tau$. Cases where  $X\leq \tau$, $\mathbb{P}-a.s.$ are of particular interest with respect to the recent literature  : in this case $Y\geq X$, $\mathbb{P}-a.s.$, and therefore, whenever both $X$ and $Y$ are $\mathbb{P}-$ integrable, it results that $$E_{\mathbb{P}}[Y| \sigma(X)] \geq X,  \ \mathbb{P}-a.s.,$$ the left hand term denoting a conditional expectation with respect to the $\sigma-$ field generated by the random variable $X$. Otherwise stated, some of those couplings are also {submartingale} couplings, while $(-X, -Y, (\Omega, \mathcal{A}, \mathbb{P}))$ is therefore a particular case of so-called supermartingale couplings, which have been recently investigated within a canonical form, for instance see \cite{NUTZ}. However, submartingale couplings don't exhaust the class of adapted couplings within the acceptation of Definition~\ref{DEFCA1}. Indeed, at the inverse, we may find random variables $X$ and $\tau$, $\mathbb{P}-$ independent from one to each other such that $X>\tau$, $\mathbb{P}-a.s.$, and take $Z=0$ in which case $X$ and $Y=\tau$ are $\mathbb{P}-$ independent, without necessarily being neither martingale nor submartingale couplings ; for a recall on definitions of martingale couplings, for instance see \cite{JUILLET}, while for other among attractive approaches using martingales together with optimal transport, for instance see \cite{BHC}. Further, as noticed earlier, Definition~\ref{DEFCA1} doesn't require $\tau$ to be necessarily independent to $X$ in $(2)$, but rather states a specific conditional independence which motivated the interpretation adopted above of $\tau$ as a waiting time for $X$.

Let $c : \mathbb{R} \times \mathbb{R} \to \mathbb{R}$ be a non-negative \textsc{Borel} measurable map, where for given $x,y\in \mathbb{R}$, $c(x,y)$ is interpreted within this model as the cost to send a signal at time $y$ when an entry signal is received at time $x$ ; for instance, when applying this model to a post office gestion, for $x<y$, the number $c(x,y)\geq 0$ may coincide with the cost, which is expressed in euros, to stock a letter or a package from the time described by $x$, where it has been received at the post office, until the time described by $y$ where it is delivered out to the target recipient. In this work, given two \textsc{Borel} probability measures $\eta$ and $\nu$ on $\mathbb{R}$, we are interested in the following problem $$\inf E_{\mathbb{P}}[c(X,Y)],$$ where the infimum is taken on the set of all the adapted couplings $(X,Y, (\Omega,\mathcal{A},\mathbb{P}))\in Cpl_a(\eta,\nu)$ of $\eta$ to $\nu$.

Here, we address this problem through the connections between those adapted couplings and so-called \textit{causal (or adapted) transport plans}. Indeed, as stated below, it turns out that adapted couplings defined according to Definition~\ref{DEFCA1} actually provide a specific representation for causal transport plans between probabilities on the set of real numbers endowed with suitable probabilistic \textit{filtrations}, within the precise definition adopted in \cite{CAUSAL}, which will be the ariane thread of this note.  This motivates the structure of this work, which is the following :

\begin{figure}
  \center
    \includegraphics[width=6cm]{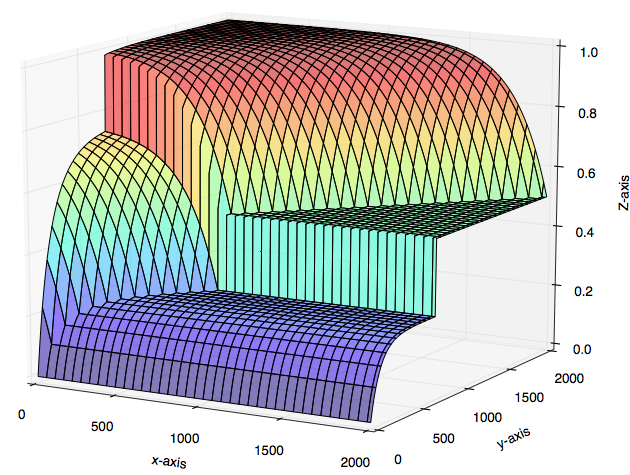}
       \includegraphics[width=6cm]{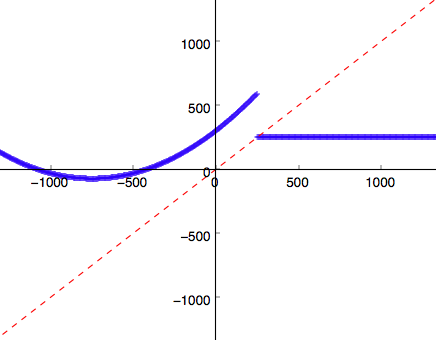}
    \caption{These figures, which will be described more accurately below, illustrate Lemma~\ref{empor} and Theorem~\ref{Theoremlk}. The image on the left represents the graph of the restriction of a conditional  cumulative distribution function $$F^\gamma :  (x,y) \in \mathbb{R}^2 \to  \mathbb{P}(\{Y \leq y\} | X = x ) \in [0,1]$$  which is associated to the joint law $p_{X,Y}$ of an  adapted coupling $(X,Y, (\Omega,\mathcal{A},\mathbb{P}))$ according to Definition~\ref{DEFCA1} : as it will be shown by Theorem~\ref{thmpm}, it turns out to coincide with some causal transport plan $\gamma$, when $\mathbb{R}$ is endowed with a suitable probabilistic filtration; above, the $"X=x"$ refers to a desintegration of measure. The image on the right represents the graph of a function $T :\mathbb{R} \to \mathbb{R}$ such that $(X,T(X), (\Omega,\mathcal{A},\mathbb{P}))$ is an adapted coupling for any real-valued random variable $X$ : as it will be pointed out below, the joint law $p_{X,T(X)}$ actually boils down to a deterministic causal transport plan. Within models, in both cases illustrated above, the upper half-plane (resp. the lower half-plane) above (resp. below) the first diagonal corresponds to cases where a response signal was sent at a time $y$ relatively \textbf{after} (resp. \textbf{before}) the time $x$ where the entry signal has been received. Depending on those two regions of the plane, curves may exhibit quite distinct aspects.}
    \end{figure}

In section \ref{section1}, we fix the notation with a particular emphasize on the filtration which we take on $\mathbb{R}$, and on the transport kernel which is used to define the filtration generated by a transport plan adopted in \cite{CAUSAL} ; within models, recall that this filtration encapsulates the \textit{information flow} which is required to perform the plan as time evolves.  In section \ref{section2}, Definition~\ref{CDefmp} provides a recall of the latter, within this specific framework, which is used to state the definition of causal transport plans which we adopt here.  Then, Lemma~\ref{empor} mentions  a characterization for such plans between \textsc{Borel} probability measures on $\mathbb{R}$, which is then used to shorten subsequent proofs : it stands on specific properties of a conditional probability cumulative distribution function $$F^\gamma : (x,y)\in \mathbb{R}^2 \to F^{\gamma}_x(y) \in \mathbb{R},$$ which is \textsc{Borel} measurable (respectively  increasing and \textit{c\`{a}d-l\`{a}g}) along the first (resp. second) variable ; the french acronym \textit{c\`{a}d-l\`{a}g} refers to right-continuous functions with left limits.  To obtain a better understanding of those latter plans within the particular circumstances which we encounter in this context, we then state an accurate  characterization to identify the deterministic ones, in Theorem~\ref{Theoremlk}, while Lemma~\ref{Lemma243} provides details to clearly grasp its origin. As a byproduct, this theorem yields a strong motivation for the representation provided by the adapted couplings adopted in Definition~\ref{DEFCA1}. A proof that this representation actually holds in the general case is provided in Theorem~\ref{thmpm} of section~\ref{section3}, where, under mild hypothesis, it is then applied to obtain, in Corollary~\ref{optimalexistence}, the existence of a solution to the problem which is stated above, from adapted tools of mass transport. It may be remarked that some arguments adopted from optimal transport can facilitate the identification of the optima. Furthermore, several examples are provided, notably from the \textit{jump times} of a simple \textsc{Poisson} process, which are aimed to provide some explicit illustrations of specificities of the related causal transport plans, within this framework.

\section{Preliminaries and notation}
\label{section1}

In this whole work, for the sake of clarity, we set $\mathbb{T}= \mathbb{R}$, while $\mathbb{R}$ is endowed with its usual topology, whose associated \textsc{Borel} $\sigma-$ field is denoted by $\mathcal{B}(\mathbb{R})$. The set of \textsc{Borel} probability measures on $(\mathbb{R}, \mathcal{B}(\mathbb{R}))$ is denoted by $M_1(\mathbb{R})$. Whenever $\mathcal{G}\subset \mathcal{B}(\mathbb{R})$ denotes a sigma-field, and $\nu\in M_1(\mathbb{R})$, the set of $\nu-$ negligible subsets is defined by $\mathcal{N}^\nu:= \{N \subset \mathbb{R} : \ \exists A\in \mathcal{B}(\mathbb{R}),  \ N\subset A, \ \text{and,} \ \nu(A)=0\}$. Moreover $\mathcal{G}^\nu$, which we still call here the $\nu-$ completion of $\mathcal{G}$,  denotes the smallest $\sigma-$ fields on $\mathbb{R}$ which contains both all the elements of $\mathcal{G}$ and all the elements of $\mathcal{N}^\nu$ ; notice that the negligible sets are taken with respect to the whole $\mathcal{B}(\mathbb{R})$, while stating the $\nu-$ completion of a sub $\sigma-$ field $\mathcal{G}$. For much details on probabilities and on specific tools as we use it subsequently, we refer to \cite{MAL1}. For $\nu\in M_1(\mathbb{R})$, the unique extension of $\nu$ to $\mathcal{B}(\mathbb{R})^\nu$ is still denoted by $\nu$, which doesn't seem to yield any confusion below.
We endow $\mathbb{R}$ with the filtration $(\mathcal{B}_t(\mathbb{R}))_{t\in\mathbb{T}}$, i.e. $\mathcal{B}_s(\mathbb{R}) \subset \mathcal{B}_t(\mathbb{R}) \subset \mathcal{B}(\mathbb{R})$, $\forall s,t\in\mathbb{T}$ such that $s\leq t$, which is defined by $$\mathcal{B}_t(\mathbb{R}) = \{A \in \mathcal{B}(\mathbb{R}) : A \subset ]-\infty, t[ \  \text{or} \ [t,+\infty[ \subset A \},$$ so that $\mathcal{B}_t(\mathbb{R})$ is equivalently defined by $\mathcal{B}_t(\mathbb{R})= \sigma(\mathcal{C}_t)$, the smallest $\sigma-$ field on $\mathbb{R}$ which contains all the subsets which constitute \begin{equation} \mathcal{C}_t = \{]-\infty, a] , a<t\}, \label{Ctdef} \end{equation} $\forall t\in \mathbb{T}$. The cartesian product $\mathbb{R}\times \mathbb{R}$ is endowed with its usual product topology which turns it into a Polish space, while the associated \textsc{Borel} sigma field is denoted by $\mathcal{B}(\mathbb{R}\times\mathbb{R})$ ; $M_1(\mathbb{R}\times \mathbb{R})$ denotes  the set of \textsc{Borel} probability measures defined on the measurable space $(\mathbb{R} \times \mathbb{R}, \mathcal{B}(\mathbb{R}\times\mathbb{R}))$, while we adapt conventional notation analogous to the case of $\mathbb{R}$ for completions.  Given $\eta,\nu\in M_1(\mathbb{R})$, in optimal transport (\cite{KANTOROVICH1}, \cite{MONGE}, \cite{PEYRE}, \cite{Vil2}), a probability $\gamma \in M_1(\mathbb{R}\times \mathbb{R})$ is said to be a \textit{transport plan} from $\eta$ to $\nu$, if and only if, its first marginal is $\eta$ (i.e. ${p_1}_\star \gamma = \eta$) and its second marginal is $\nu$ (i.e. ${p_2}_\star\gamma =\nu$), $p_1 : (x,y)\in \mathbb{R}\times \mathbb{R} \to x\in \mathbb{R}$ (resp. $p_2 : (x,y)\in \mathbb{R}\times \mathbb{R} \to y\in \mathbb{R}$) denoting the canonical projections functions, while ${p_1}_\star\gamma$ (resp. ${p_2}_\star \gamma$) denotes the \textit{pushforward} (see \cite{MAL1}) of $\gamma$ by $p_1$ (resp. by $p_2$) ; ${p_1}_\star\gamma(A) = \gamma(A\times \mathbb{R})$ (resp.  ${p_2}_\star\gamma(B) = \gamma( \mathbb{R}\times B)$), $\forall A\in \mathcal{B}(\mathbb{R})$ (resp. $\forall B\in \mathcal{B}(\mathbb{R})$). Then,  $$\Pi(\eta, \nu) = \{\gamma \in M_1(\mathbb{R}\times \mathbb{R}) : {p_1}_\star \gamma = \eta, {p_2}_\star \gamma = \nu \}$$ denotes the set of \textit{transport plans} from $\eta$ to $\nu$.

 To describe accurately cases where the mass is splitted within the transport, recall that for any $\gamma\in \Pi(\eta,\nu)$, there exists a function $x\in \mathbb{R} \to Q_{\gamma}^x\in M_1(\mathbb{R})$ (see \cite{DM3}, \cite{I-W}), such that the function 
 \begin{equation} \label{phibdef} \phi_B : x\in \mathbb{R} \to Q^x_{\gamma}(B) \in \mathbb{R} \end{equation} is \textsc{Borel} measurable, $\forall B\in \mathcal{B}(\mathbb{R})$, 
 and   \begin{equation} \gamma(A\times B) = \int_A Q_\gamma^x(B) \eta(dx), \ \forall A,B\in\mathcal{B}(\mathbb{R}), \label{eqmle4}  \end{equation} the latter denoting a \textsc{Lebesgue} integral with respect to the \textsc{Borel} probability measure $\eta$, which is well defined from the previous hypothesis on~(\ref{phibdef}) ; within mass transport models, $Q^x_\gamma(B)$ may be interpreted as the proportion of the mass located at $x\in\mathbb{R}$ which is carried to $B\in\mathcal{B}(\mathbb{R})$, since $Q_\gamma$ coincides with a desintegration kernel of $\gamma$ with respect to the canonical projection  $p_1$ (see \cite{DM}), outside an $\eta-$ negligible set ; we have, for any $ A\in \mathcal{B}(\mathbb{R})$, $Q_\gamma^x(A)= \gamma(\{p_2 \in A\} | p_1=x )$, $\eta-a.s.$, where $"p_1=x"$ refers here to a desintegration of measure. Subsequently, any such function $Q_\gamma$ will be called a transport kernel associated to $\gamma$ ; it is unique up to an $\eta-$ negligible set, which will justify to take $\eta-$ completions below.  Futher, for $x\in \mathbb{R}$, $F^\gamma_x$ denotes the cumulative distribution function of $Q_\gamma^x$ ; $F^\gamma_x(a)= Q_\gamma^x(]-\infty,a])$, $\forall a\in \mathbb{R}$, and we set  $$F^\gamma : (x,y)\in \mathbb{R}^2\to F^\gamma_x(y)\in \mathbb{R},$$ to which we refer to as \textit{a} conditional cumulative distribution function associated to $\gamma$, as it may depend on the choice of $Q_\gamma$. Furthermore, given a real-valued random variable $X$ which is defined on a probability space $(\Omega, \mathcal{A}, \mathbb{P})$, that is, a function $X: \Omega \to \mathbb{R}$ whose level sets $\{\omega\in \Omega : X(\omega) \leq  l \} \in \mathcal{A}$ are all contained in $\mathcal{A}$, $\forall l\in \mathbb{R}$, for short and once given $A\subset \mathbb{R}$, we may use the notation $\{X\in A\}$ to denote the inverse image $X^{-1}(A):=\{\omega\in \Omega : X(\omega) \in A\}$, while we use the notation $E_{\mathbb{P}}[X]$ to denote the mathematical expectation of $X$ under the probability $\mathbb{P}$ ; $E_{\mathbb{P}}[X] := \int_{\Omega} X(\omega) \mathbb{P}(d\omega)$ denotes a \textsc{Lebesgue} integral with respect to the probability measure $\mathbb{P}$ whenever $X$ is $\mathbb{P}-$ integrable (i.e. $E_{\mathbb{P}}[|X|]<+\infty$), in which case $E_{\mathbb{P}}[X]\in \mathbb{R}$, while we set  $E_{\mathbb{P}}[X]= +\infty$ whenever $X\geq 0$, $\mathbb{P}-a.s.$ and $X$ is not $\mathbb{P}-$ integrable. The law of the random variable $X$ is denoted by $p_X\in M_1(\mathbb{R})$, where $p_X(A)= \mathbb{P}(X^{-1}(A))= \mathbb{P}(\{X\in A\})$, $\forall A\in\mathcal{B}(\mathbb{R})$, while further assuming that $Y :\Omega\to \mathbb{R}$ denotes a random variable defined on the same probability space as $X$, $p_{X,Y} \in M_1(\mathbb{R}\times \mathbb{R})$ denotes the joint law of $X$ and $Y$, that is, the pushforward $p_{X,Y}:=(X,Y)_\star \mathbb{P}$ of the probability measure $\mathbb{P}$ by the measurable function $(X,Y) : \omega\in \Omega \to (X(\omega), Y(\omega)) \in \mathbb{R}\times \mathbb{R}$, so that
 $p_{X,Y}(A)= \mathbb{P}(\{(X,Y)\in A\})$, $\forall A\in \mathcal{B}(\mathbb{R}\times \mathbb{R})$; whenever $X$ is $\mathbb{P}-$ independent to $Y$ we denote by $p_X\otimes p_Y$ their joint law $p_{X,Y}$. Finally, given a probability space $(\Omega, \mathcal{A}, \mathbb{P})$, if $C\in\mathcal{A}$ is such that $\mathbb{P}(C)>0$, we apply the standard notation $\mathbb{P}(A|C) = \frac{\mathbb{P}(A\cap C)}{\mathbb{P}(C)}$, $\forall A \in\mathcal{A}$, as it is used in Definition~\ref{DEFCA1} with $C=\{X\geq t\}$, $t\in\mathbb{T}$, while for any $A\in\mathcal{A}$, $\mathds{1}_A : \Omega \to \mathbb{R}$ denotes here the indicator function of $A$ ; $\mathds{1}_A(\omega) =1$ if $\omega \in A$, while $\mathds{1}_A(\omega) =0 $ if $\omega \notin A$, $\forall \omega \in \Omega$.

\section{Causal mass transport between \textsc{Borel} probabilities on the set of real numbers.} 

\label{section2}

\subsection{The filtration generated by a transport plan, a characterization}

We recall some conventional terminology adopted in \cite{CAUSAL}, which will be required subsequently. 

\begin{definition} \label{CDefmp}
Let $\eta,\nu\in M_1(\mathbb{R})$, and $\gamma\in \Pi(\eta,\nu)$. The filtration $(\mathcal{G}_t(\gamma))_{t\in \mathbb{T}}$ generated by $\gamma$ is defined, for any $t\in\mathbb{T}$,  by $$\mathcal{G}_t(\gamma) = \sigma(\phi_B : B\in \mathcal{B}_t(\mathbb{R}), \ \nu(\partial B)=0)^\eta ;$$ the $\eta-$ completion of the smallest sigma field which turns into measurable functions all the functions $\phi_B$ as in~(\ref{phibdef})  such that $B\in \mathcal{B}_t(\mathbb{R})$ and $\nu(\partial B)=0$. The set $\Pi_c(\eta,\nu)$ of causal (or adapted) transport plans $\gamma$  from $\eta$ to $\nu$ is defined by  $$\Pi_c(\eta,\nu) = \{\gamma \in \Pi(\eta,\nu) :  \mathcal{G}_t(\gamma) \subset \mathcal{B}_t(\mathbb{R})^\eta, \ \forall t\in \mathbb{T} \}.$$
\end{definition}

\begin{lemma} \label{empor}
Let $\eta,\nu\in M_1(\mathbb{R})$,  $\gamma\in \Pi(\eta,\nu)$, and further denote by $Q_\gamma : x \in \mathbb{R} \to Q_\gamma^x \in  M_1(\mathbb{R})$ a transport kernel associated to $\gamma$ by~(\ref{eqmle4}). Then the following assertions are equivalent :
\begin{enumerate}[(i)]
\item $\gamma \in \Pi_c(\eta,\nu)$
\item For any $t\in \mathbb{T}$ such that $\eta([t,+\infty[)>0$, $\forall A\in \mathcal{B}_t(\mathbb{R})$,and  for all $x\geq t$ outside an $\eta-$ negligible set, we have $$Q_\gamma^x(A) = \gamma (\{p_2 \in A\} | \{p_1 \in [t,+\infty[\});$$ conditional probability, with respect to an event.
\item For any $t\in \mathbb{T}$  such that $\eta([t,+\infty[)>0$, and for any $a<t$, $\exists c_{a,t}\in \mathbb{R}$ such that we have $$F_x^\gamma(a) = c_{a,t}, \ \forall x \geq t, \ \text{outside  an} \ \eta-\text{negligible set} ,$$ $F_x^\gamma : \mathbb{R} \to [0,1]$ denoting the cumulative distribution function of $Q_\gamma^x$, $\forall x\in \mathbb{R}$.
\end{enumerate}
\end{lemma}
\nproof
Since the definitions yield $$\mathcal{B}_t(\mathbb{R}) =\rho_t^{-1}(\mathcal{B}(\mathbb{R})) =\{ \{x\in \mathbb{R} : \rho_t(x)\in A \} : A\in \mathcal{B}(\mathbb{R}) \},$$ where the function $\rho_t : s\in \mathbb{R} \to \min(s,t)\in \mathbb{R}$ is continuous, $\forall t\in \mathbb{T}$, and $\rho_t\circ \rho_u = \rho_{\min(t,u)},$ $\forall t,u\in \mathbb{T}$, $\circ$ denoting the usual pullback of functions, from \cite{CAUSAL} we first deduce that, for any $t\in \mathbb{T}$, we have \begin{equation}\mathcal{G}_t(\gamma) = \sigma(\phi_A : A\in \mathcal{B}_t(\mathbb{R}))^\eta ; \label{eqfiltregal}  \end{equation}  since $\mathcal{B}_t(\mathbb{R})$ is also the $\sigma-$ field generated by $\{\rho_t^{-1}(F) : \ F \ \text{closed  subset of} \  \mathbb{R} \}$, the latter equation~(\ref{eqfiltregal}) can also be  checked readily by a monotone class argument, once noticed that, since $\rho_t$ is continuous, for any closed set $F\subset\mathbb{R}$,  the holding condition $\lambda(\{y>0 : \nu(\partial (\rho_t^{-1}(d_F^{-1}([0,y]))))>0\})=0$,  where $\lambda$ denotes the  \textsc{Lebesgue} mesure, and where $d_F : x \in \mathbb{R} \to \inf_{y\in F} |x-y| \in \mathbb{R}$, yields the existence of a decreasing sequence $ (B_n)_{n\in\mathbb{N}}$  of closed subsets of the form $\rho_t^{-1}(d_F^{-1}(]-\infty,t_n]))$, for some $t_n>0$, which further satisfy $\nu(\partial B_n)=0$, $B_n\in \mathcal{B}_t(\mathbb{R})$, $\forall n\in \mathbb{N}$ while $\rho_t^{-1}(F) =\bigcap_{n\in\mathbb{N}} B_n$,  so that  $\phi_{\rho_t^{-1}(F)}=\lim_{n\to \infty} \phi_{B_n}$, $\eta-a.s.$, and therefore $\phi_{\rho_t^{-1}(F)}$ is $\mathcal{G}_t(\gamma)-$ measurable.  From~(\ref{eqfiltregal}), we first notice that $\gamma\in \Pi_c(\eta,\nu)$, if an only if, \begin{equation} Q_\gamma(A) = E_\eta[ Q_\gamma(A) | \mathcal{B}_t(\mathbb{R})^\eta], \ \eta-a.s. \label{qgamma1e} \end{equation} for all $A\in\mathcal{B}_t(\mathbb{R})$, $\forall t\in\mathbb{T}$ ; the right-hand term of~(\ref{qgamma1e}) denotes a conditional expectation with respect to the $\sigma$- field $\mathcal{B}_t(\mathbb{R})^\eta$, where $t\in\mathbb{T}$.  We now turn to the main part of the proof. First assume $(i)$, and let  $t\in\mathbb{R}$. For $A\in \mathcal{B}_t(\mathbb{R})$, from the very definition of the conditional expectation with respect to a sigma field (for instance see \cite{DM}, \cite{I-W}), and from the properties of the filtration $(\mathcal{B}_t(\mathbb{R}))_{t\in\mathbb{T}}$ just recalled above,  it directly follows that  
$$E_\eta[ Q_\gamma(A) | \mathcal{B}_t(\mathbb{R})^\eta] = Q_\gamma(A) \mathds{1}_{]-\infty,t[} + \gamma(\{p_2 \in A\} | \{p_1\in[t,+\infty[\}) \mathds{1}_{[t,+\infty[} , \ \eta-a.s.,$$ if $\eta([t,+\infty[)>0$, while~(\ref{qgamma1e}) already holds as soon as $\eta([t,+\infty[)=0$. As a consequence,~(\ref{qgamma1e}) holds for all $A\in\mathcal{B}_t(\mathbb{R})$, and $\forall t\in\mathbb{T}$, if and only if,  $$ \gamma(\{p_2 \in A\} | \{p_1\in[t,+\infty[\})\mathds{1}_{[t,+\infty[} = Q_\gamma(A) \mathds{1}_{[t,+\infty[}, \ \eta-a.s.,$$ holds,  for all $A\in\mathcal{B}_t(\mathbb{R})$, and $\forall t\in\mathbb{T}$ such that $\eta([t,+\infty[)>0$, which is $(ii)$ : we have proved that $(i)$ is equivalent to $(ii)$. On the other hand, since for $t\in \mathbb{T}$ and $a<t$, $]-\infty, a] \in\mathcal{B}_t(\mathbb{R})$, $(ii)$ yields $(iii)$ as a particular case, by taking \begin{equation} c_{a,t}= \gamma(\{p_2 \in ]-\infty,a]\} | \{p_1\in[t,+\infty[\})\in [0,1].  \label{catformula} \end{equation} Finally, assuming that $(iii)$ holds, for  $t\in \mathbb{T}$  such that $\eta([t,+\infty[)>0$, and $a<t$, we notice that we still obtain~(\ref{catformula}), so that we have $(ii)$ for all $A\in \mathcal{B}_t(\mathbb{R})$ of the form $]-\infty,a]$, $a<t$. Since $\mathcal{B}_t(\mathbb{R})= \sigma(\mathcal{C}_t)$, where $\mathcal{C}_t\subset\mathcal{B}(\mathbb{R})$ is given by~(\ref{Ctdef}),  $(ii)$ follows from a monotone class argument (see \cite{DM}).   \nqed

\begin{figure}
   \center
    \includegraphics[width=10cm]{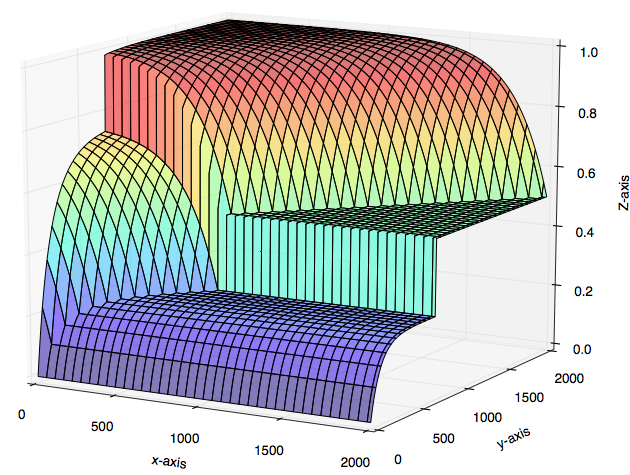}
    \caption{ This figure illustrates Lemma~\ref{empor}. Representation of the restriction to the square $[0 ,2000]\times [0,2000]$ of a conditional  cumulative distribution function $F^\gamma$, which is associated to the causal transport plan  $$\gamma =
      \frac{1}{2} p_{X,X+Z} + \frac{1}{4} p_X\otimes \delta_{t_0} + \frac{1}{4}p_X\otimes p_X,$$ where $t_0=700$, and where $X$, and $Z$ denote two $\mathbb{P}-$ independent real valued random variables defined on a same probability space $(\Omega, \mathcal{A}, \mathbb{P})$, whose laws are given by  $p_X= \mathcal{E}\left(\frac{1}{100}\right)$, and by  $p_Z =\mathcal{E}\left(\frac{1}{200}\right)$ ; $\mathcal{E}(\beta)\in M_1(\mathbb{R})$ denotes the so-called \textit{exponential law} of parameter $\beta\in \mathbb{R}_+^\star$, that is, the \textsc{Borel} probability measure absolutely continuous with respect to the \textsc{Lebesgue} measure, whose density $f_\beta : \mathbb{R} \to \mathbb{R}$ is given by $$f_\beta(x)= \beta  \exp(-\beta x) \mathds{1}_{[0,+\infty[}(x), \ \forall x\in \mathbb{R}.$$ }
\end{figure}

\subsection{Deterministic causal transport plans}

Let $\eta\in M_1(\mathbb{R})$, let $T : \mathbb{R} \to \mathbb{R}$ be a \textsc{Borel} measurable function, and define $\gamma_T = (I_\mathbb{R}, T)_{\star}\eta$,  the deterministic transport plan of $\eta$ by $T$ ; we can set $Q^x_{\gamma_T} =\delta_{T(x)}$, $\forall x\in \mathbb{R}$, $\delta_y$ denoting the \textsc{Dirac} mass (see \cite{MAL1}) centered on $y\in \mathbb{R}$, $\forall y\in \mathbb{R}$,  while  $\gamma_T$ coincides with the pushforward measure of $\eta$ by the measurable function $$(I_{\mathbb{R}}, T) : x\in \mathbb{R} \to (x,T(x))\in \mathbb{R}\times \mathbb{R},$$ so that $(I_{\mathbb{R}}, T)_\star {\eta} (C) = \eta(\{x\in \mathbb{R} : (x,T(x))\in C\})$, $\forall C \in \mathcal{B}(\mathbb{R}\times\mathbb{R})$. Further recall that, within those specific hypothesis,  from~(\ref{eqfiltregal})  we obtain \begin{equation} \label{GTDEF}  \mathcal{G}_t(\gamma_T) = \mathcal{G}_t^T,\end{equation} where   $$\mathcal{G}_t^T = \{T^{-1}(A) \ :\ A\in \mathcal{B}_t(\mathbb{R})\}^\eta,$$ $\forall t\in \mathbb{T}$  ; $T^{-1}(A) = \{x \in \mathbb{R} : T(x)\in A\}$ still denotes the inverse image of $A\in\mathcal{B}(\mathbb{R})$ by the function $T$. Finally, for the sake of consistency with \cite{CAUSAL} and for short, given $\eta\in M_1(\mathbb{R})$, subsequently $T$ will be said to be \textit{adapted} (or \textit{causal}), if and only if, $\gamma_T\in \Pi_c(\eta, \nu)$, where $\nu=T_\star \eta$, if and only if,  $$\mathcal{G}_t^T\subset \mathcal{B}_t(\mathbb{R})^\eta, \ \forall t\in \mathbb{T},$$ where the last equivalence follows from~(\ref{GTDEF}); again, $T_\star\eta\in M_1(\mathbb{R})$ denotes the pushforward or direct image of $\eta\in M_1(\mathbb{R})$ by the \textsc{Borel} measurable function $T : \mathbb{R} \to \mathbb{R}$, that is $T_\star \eta(A)= \eta(\{x\in\mathbb{R} : T(x)\in A\})$, $\forall A\in\mathcal{B}(\mathbb{R})$. 

\begin{lemma} \label{Lemma243}
Let $\eta\in M_1(\mathbb{R})$, and let $T : \mathbb{R} \to \mathbb{R}$ be a \textsc{Borel} measurable function such  that $T$ is adapted, i.e. $T^{-1}(\mathcal{B}_t(\mathbb{R})) \subset \mathcal{B}_t(\mathbb{R})^\eta$, $\forall t\in \mathbb{T}$. Further assume that we have  \begin{equation} \eta(\{x\in \mathbb{R} : T(x)< x\}) >0 . \label{hyptinf1} \end{equation} Then, the following assertions are satisfied :

\begin{enumerate}[(i)]
\item $\forall t\in \mathbb{T}$, $\forall A\in \mathcal{B}(\mathbb{R})$ such that either $A\subset ]-\infty,t[$ or $[t,+\infty[\subset A$, the following implication holds : $$ \eta(\{x \geq t : T(x) \in A \}) >0 \implies \eta(\{x\in \mathbb{R} : T(x) \in A\} | [t,+\infty[) =1.$$
\item Let $$D = \{t \in \mathbb{R} : \eta(\{ x \geq t  :  T(x) <t \}) >0\}.$$ Then, $D\neq \emptyset$,  and there exists a unique $t_0\in \mathbb{R}$ which is such that $t_0 < t$, $\forall t\in D$,  and  $$\eta(\{x\in \mathbb{R} : T(x) =t_0\} |[t,+\infty[)= 1, \ \forall t\in D.$$ 
\item  $\inf D \in \mathbb{R}$ is not attained and we have $$\eta(]\inf D, +\infty[) >0$$ and  $$\eta(\{x \in \mathbb{R} : T(x) =t_0 \} | ]\inf D,+\infty[) =1.$$
\item  The following hold :  \begin{equation*} \begin{cases}  \eta(\{x \leq t_0 : T(x)\geq x\}) = \eta(]-\infty, t_0])  \\  \eta(]t_0,+\infty[) >0 
 \\ \eta(\{x\in \mathbb{R} : T(x) =t_0\} | ]t_0,+\infty[) =1 
 \end{cases} . \end{equation*} 
\end{enumerate}
\end{lemma}
\nproof
Since $\gamma_T =(I_\mathbb{R}, T)_\star \eta \in \Pi_c(\eta,\nu)$, where $\nu=T_\star \eta$, Lemma~\ref{empor} applies. Hence, assume that $t\in \mathbb{T}$ and $A\in \mathcal{B}_t(\mathbb{R})$ are such that $\eta(\{x\geq t : T(x)\in A\})>0$ ; in particular, $\eta([t,+\infty[)>0$. As we can state $$Q^x_{\gamma_T}(A) = \delta_{T(x)}(A) =\mathds{1}_A(T(x)) = \mathds{1}_{T^{-1}(A)}(x), \ \forall x\in \mathbb{R}, $$ $Q_{\gamma_T}$ denoting a transport kernel for $\gamma_T$ recalled above, from Lemma~\ref{empor} we obtain $$\mathds{1}_{T^{-1}(A)} \mathds{1}_{[t,+\infty[} = \gamma_T( \{p_2\in A\} | \{p_1\in[t,+\infty[\}) \mathds{1}_{[t,+\infty[}, \ \eta-a.s.,$$ so that   $\eta(\{x \geq t : T(x) \in A \}) >0$  implies $ \gamma_T( \{p_2\in A\} | \{p_1\in[t,+\infty[ \}) =1$, from which the implication $(i)$ follows. Further noticing that  $$\eta(\{x\in \mathbb{R} : T(x)< x\}) = \eta\left(\bigcup_{q\in\mathbb{Q}}\{ x\geq q : T(x) < q \} \right),$$ and that a countable union of negligible sets is again negligible, we obtain $D\neq \emptyset$ from~(\ref{hyptinf1}). Then, we choose any $t\in D$, and define $$t_0= \inf(\{ y\in \mathbb{R} : \eta(\{x\geq t : T(x) <y \}) \geq \eta(\{x\geq t : T(x) <t\}) \}).$$ Together with $(i)$, from the definitions, it is then an easy exercise to check that $t_0\in\mathbb{R}$ doesn't depend on the choice of $t\in D$, while it meets the required hypothesis.  Finally, the rest of the statement results by  following the successive steps, while applying $(i)$ as far as necessary.
 \nqed

\begin{theorem} \label{Theoremlk}
Let $\eta\in M_1(\mathbb{R})$, let $T : \mathbb{R} \to \mathbb{R}$ be a \textsc{Borel} measurable function, set $\nu = T_\star \mu$, and still define $\gamma_T = (I_\mathbb{R}, T)_{\star}\eta$ to be the deterministic transport plan of $\eta$ by $T$. Then, $\gamma_T \in \Pi_c(\eta,\nu)$, if and only if, $T$ satisfies one among the two hypothesis below :
\begin{enumerate}
\item $T(x)\geq x$, $\eta-a.s.$
\item $ \exists t_0 \in \mathbb{R}$ which is such that, $\forall x \leq t_0$ outside an $\eta-$ negligible set, we have $T(x)\geq x$, and for any $x>t_0$ outside an $\eta-$ negligible set, we have $T(x)= t_0$.

\end{enumerate}

\end{theorem}
\nproof
Assuming that one among  $(1)$ or $(2)$ is satisfied, it is enough to notice that $\mathcal{B}_t(\mathbb{R}) =\sigma(\mathcal{C}_t)$, where $\mathcal{C}_t$ is given by~(\ref{Ctdef}), so that we easily check  that $\gamma_T \in \Pi_c(\eta,\nu)$, as $\mathcal{G}_t^T \subset \mathcal{B}_t(\mathbb{R})^\eta$, $\forall t\in \mathbb{T}$ ; recall that~(\ref{GTDEF}) holds. Conversely, assuming that $\gamma_T \in \Pi_c(\eta,\nu)$, the result follows from Lemma~\ref{Lemma243}. 
\nqed

\begin{remark} $(2)$ above is clearly equivalent to : 
$\exists t_0\in \mathbb{R}$ and a non-negative \textsc{Borel} measurable function $g_T : \mathbb{R} \to [0,+\infty[$, such that $$T(x) = (x+ g_T(x))\mathds{1}_{]-\infty,t_0]}(x)  + t_0 \mathds{1}_{]t_0,+\infty[}(x), \ \eta-a.e..$$
\end{remark}

\begin{figure}
   \center
    \includegraphics[width=3cm]{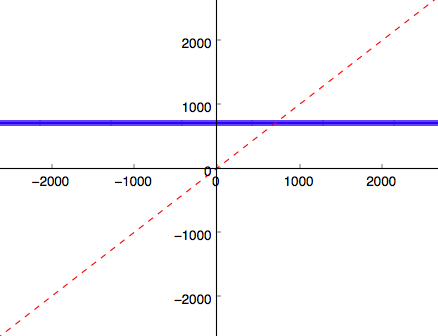}
    \includegraphics[width=3cm]{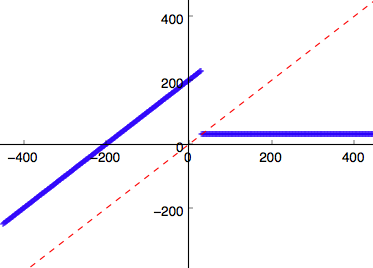}
     \includegraphics[width=3cm]{AOTFig3}
     \includegraphics[width=3cm]{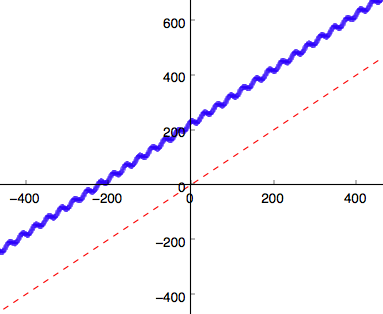}
    \caption{This figures illustrate Theorem~\ref{Theoremlk}. In red, the first diagonal of the plane. In blue, the representation of the restriction of the graph of the \textsc{Borel} measurable function $T : \mathbb{R} \to \mathbb{R}$, given by, from the image on the left to the image on the right,
$$T(x) = c, \ \forall x\in \mathbb{R},$$ 
with $c=700$, by  $$T(x) = (x+z_0)\mathds{1}_{]-\infty,t_0]}(x)+ t_0 \mathds{1}_{]t_0,+\infty[}(x), \ \forall x\in \mathbb{R},$$ with $t_0=30$ and $z_0=200$, by $$T : x\in \mathbb{R} \to (ax^2+bx +c) \mathds{1}_{]-\infty,t_0]} + t_0 \mathds{1}_{]t_0,+\infty[}(x) \in \mathbb{R},$$ with $t_0=250$, $a=\frac{1}{1500}$, $b=1$, and $c=300$, and finally by $$T : x\in \mathbb{R} \to  x+ a \cos( b x) + c\in \mathbb{R},$$ where $a=10$, $b= \frac{1}{5}$, and $c=220$. For any $\eta\in M_1(\mathbb{R})$, those functions $T$ induce a deterministic causal transport plan $\gamma_T =(I_{\mathbb{R}}, T)_\star \eta\in \Pi_c(\eta,\nu)$ of $\eta$ to $\nu=T_\star\eta$. } 
    
    \end{figure}

As a first practical consequence of Theorem~\ref{Theoremlk}, we get  :

\begin{corollary}
Let $T :\mathbb{R} \to \mathbb{R}$ be a \textsc{Borel} measurable function, $\eta\in M_1(\mathbb{R})$, and assume that $\nu=T_\star\eta$ and $\eta$ are diffuse (i.e. $\eta(\{x\})= \nu(\{x\})=0$, $\forall x\in \mathbb{R}$), and  that $\eta(]t,+\infty[)>0$, $\forall t\in \mathbb{T}$ ; for instance, $\eta$ has a cumulative distribution function which is strictly increasing on $]y,+\infty[$ for some $y\in \mathbb{R}$. Then, $T$ is adapted, if and only if, $$T(x)\geq x, \ \eta-a.s..$$  \nqed
\end{corollary}

\begin{example}
Let $a,b\in \mathbb{R}$, and define $T : x\in \mathbb{R} \to a x + b\in \mathbb{R}$ ; recall again that within those hypothesis, given $\eta\in M_1(\mathbb{R})$,  $T$ is said to be adapted, if and only if, $\gamma = (I_\mathbb{R}, T)_{\star}\eta\in \Pi_c(\eta,\nu)$, where $\nu=T_\star\eta$, if and only if, $T^{-1}(\mathcal{B}_t(\mathbb{R}))\subset \mathcal{B}_t(\mathbb{R})^\eta$, $\forall t\in \mathbb{T}$.  
\begin{itemize}
\item Assuming that $\eta = \mathcal{N}(0,1)$ (the \textsc{Gauss} law), $T$ is adapted, if and only if, either $a=0$, or both $a=1$ and $b\geq 0$.
\item If $\eta = \mathcal{E}(1)$ (\textit{exponential law}) then $T$ is adapted, if and only if,  either $a=0$, or both $a\geq1$ and $b\geq 0$.
\end{itemize}
\end{example}

\section{Characterization of adapted couplings and probabilistic optimization.}
\label{section3}

The following Theorem~\ref{thmpm} will allow to address the problem stated in the introduction by using tools of causal optimal transport, in Corollary~\ref{optimalexistence}.

\begin{theorem} \label{thmpm}  
\begin{enumerate}
\item Let $\eta,\nu \in M_1(\mathbb{R})$. Then, for any $\gamma\in \Pi_c(\eta, \nu)$, there exists an adapted coupling $(X, Y, (\Omega, \mathcal{A},\mathbb{P}))\in Cpl_a(\eta,\nu)$ such that $p_{X,Y}= \gamma$, where $p_{X,Y}\in M_1(\mathbb{R}\times \mathbb{R})$ denotes the joint law of the random variables $X :\Omega \to \mathbb{R}$ and $Y : \Omega \to \mathbb{R}$ on the probability space $(\Omega,\mathcal{A}, \mathbb{P})$, and where $Cpl_a(\eta,\nu)$ (resp. $\Pi_c(\eta,\nu)$) is given by Definition~\ref{DEFCA1} (resp. by Definition~\ref{CDefmp}). 
\item For any adapted coupling $(X, Y, (\Omega, \mathcal{A},\mathbb{P}))\in Cpl_a(p_X,p_Y)$, we have $$p_{X,Y}\in \Pi_c(p_X,p_Y),$$ $p_{X}$ (resp. $p_Y$) denoting the probability law of $X$ (resp. of $Y$) on $(\Omega,\mathcal{A}, \mathbb{P})$.
\end{enumerate}
\end{theorem}
\nproof
First assume that $\gamma\in \Pi_c(\eta,\nu)$, and define \begin{equation*} \begin{cases}  (\Omega, \mathcal{A}, \mathbb{P}) =  (\mathbb{R}\times \mathbb{R}, \mathcal{B}(\mathbb{R}\times \mathbb{R})^\gamma, \gamma)  \\  X   =  p_1 \\ Y  =  p_2 \\ \tau  = \min(p_1,p_2) 
 \\ Z  =  (Y-X) \mathds{1}_{\{X\leq \tau\}} 
 \end{cases} ; \end{equation*} 
 $p_1$ (resp. $p_2$) still denotes the canonical projection on the first (resp. second) component of the product space $\mathbb{R}\times\mathbb{R}$. Then, it is enough to check that $(X,Y, (\Omega,\mathcal{A},\mathbb{P}))$ satisfies the axioms of Definition~\ref{DEFCA1} to obtain $(1)$, which follows from Lemma~\ref{empor}. Then,  $(2)$ follows by applying again Lemma~\ref{empor}, $(iii)$, together with the definitions. 
 \nqed

\begin{figure}
  \center
    \includegraphics[width=6cm]{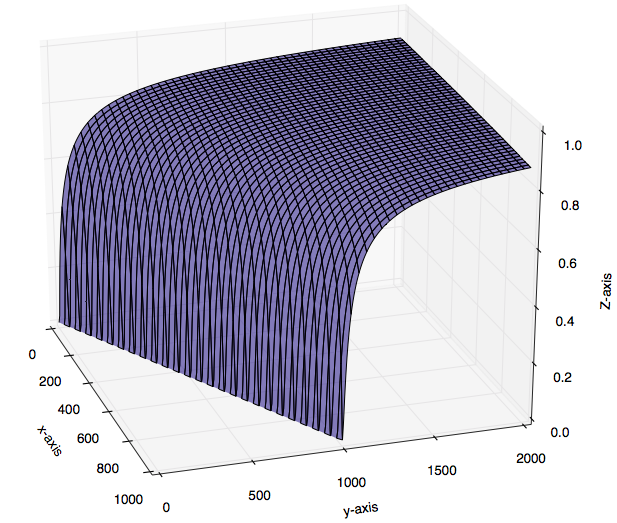}
  \caption{This figure illustrates Example~\ref{brownianexemple}, for $0<a<b$.  From classical results (see 2/ in 1.7. of \cite{IM}), for $\gamma = p_{T_a,T_b}$ it follows that we can choose $$F^\gamma_x(y) = \mathds{1}_{]0,+\infty[}(y-x)\left(1 - \text{erf}\left(\frac{b-a}{\sqrt{2(y-x)}}\right) \right),$$ $\forall (x,y)\in \mathbb{R}_+\times \mathbb{R}$, where \textit{erf} denotes the so-called error function, which is given by $\text{erf}(z) = \frac{2}{\sqrt{\pi}}\int_0^z \exp\left(-t^2\right) dt$, $\forall z\in \mathbb{R}_+$; $\mathbb{R}_+$ is the support of $p_{T_a}$. The figure represents the restriction of the graph of $F^\gamma$, for $a=6$ and $b=11$.}
    \end{figure}
\begin{example} \label{brownianexemple}
Let $(B_t)_{t\in[0,+\infty[}$ be a standard real valued Brownian motion, starting from $0$, which is defined on a $\mathbb{P}-$ complete probability space $(\Omega, \mathcal{A}, \mathbb{P})$ (see \cite{PROTTER});  this is the continuous \textsc{L\'{e}vy} process associated to the \textsc{Gauss} law $\mathcal{N}(0,1)$ which is an infinitely divisible distribution, while the law of this continuous stochastic process is the so-called \textsc{Wiener} measure, see  \cite{MALLIAVIN}, \cite{PROTTER}, \cite{STROOCK}, \cite{WIENER1}. In particular, $(B_t)_{t\in [0,+\infty[}$ is a family of gaussian random variables labeled by $\mathbb{R}_+$, such that on the one hand the function $t \in [0,+\infty[ \to B_t(\omega)\in \mathbb{R}$ is continuous for any $\omega \in \Omega$ outside a $\mathbb{P}-$ negligible set, while on the other hand, for any $s,t\in ]0,+\infty[$ such that $s<t$,  there exists two independent random variables $G$ and $\widetilde{G}$, depending respectively on $s$ and on both $s$ and $t$, defined on a same probability space, and with a same \textsc{Gauss} law $p_G=p_{\widetilde{G}}=\mathcal{N}(0,1)\in M_1(\mathbb{R})$, such that the joint law $p_{B_s,B_t}$ of the two random variables $B_s :\Omega \to \mathbb{R}$ and $B_t : \Omega \to \mathbb{R}$ on $(\Omega, \mathcal{A}, \mathbb{P})$ coincide with the joint law of the random variables $\sqrt{s} G$ and $\sqrt{s}G + \sqrt{t-s} \widetilde{G}$, that is $$p_{B_s, B_t}= p_{\sqrt{s} G, \sqrt{s}G + \sqrt{t-s} \widetilde{G}}.$$  For $c\in \mathbb{R}$, define the, $\mathbb{P}-$ almost surely finite, first passage time random variable  $$T_c := \inf(\{ t \geq 0 : B_t = c \}),$$ which is a so-called stopping time with respect to the filtration generated by $(B_t)_{t\in[0,+\infty[}$ (see \cite{IM}).  Assuming that $0<a<b$, from Theorem~\ref{thmpm}, we obtain $$p_{T_a,T_b}\in \Pi_c(p_{T_a}, p_{T_b}),$$ $p_{T_a}$ (resp. $p_{T_b}$) denoting the probability law of $T_a$ (resp. of $T_b$), while $p_{T_a,T_b}= (T_a, T_b)_\star \mathbb{P} \in M_1(\mathbb{R}\times \mathbb{R})$ denotes their joint law, on the probability space $(\Omega, \mathcal{A}, \mathbb{P})$. 
\end{example}

\begin{corollary}
\label{optimalexistence}
Assume that $c : \mathbb{R}\times \mathbb{R} \to \mathbb{R}$ is a non-negative and lower semi-continuous function, and let $\eta,\nu \in M_1(\mathbb{R})$. Then, there exists an adapted coupling $(X_\star,Y_\star, (\Omega_\star,$ $\mathcal{A}_\star,$ $\mathbb{P}_\star))\in Cpl_a(\eta,\nu)$ which attains \begin{equation} \label{AC23plpb1} \inf(\{E_{\mathbb{P}}[c(X,Y)] : (X,Y ,(\Omega, \mathcal{A}, \mathbb{P})) \in Cpl_a(\eta,\nu)\}), \end{equation} where $Cpl_a(\eta,\nu)$ denotes the set defined in Definition~\ref{DEFCA1}.  Moreover, this problem is equivalent to the causal (or adapted) \textsc{Monge}-\textsc{Kantorovich} problem \begin{equation} \inf_{\gamma\in \Pi_c(\eta,\nu)} \int_{\mathbb{R}\times\mathbb{R}} c(x,y) \gamma(dx,dy), \label{MKC1P} \end{equation} where $\Pi_c(\eta,\nu)$ is given by Definition~\ref{CDefmp}. That is, $\gamma \in \Pi_c(\eta,\nu)$ attains the infimum of~(\ref{MKC1P}), if and only if, there exists some adapted coupling $(X,Y, (\Omega,\mathcal{A},\mathbb{P}))\in Cpl_a(\eta,\nu)$ which attains the infimum of~(\ref{AC23plpb1}), and $\gamma =p_{X,Y}$,  where $p_{X,Y}$ denotes the joint law of the pair of random variables $X : \Omega \to \mathbb{R}$ and $Y: \Omega \to \mathbb{R}$, on the $\mathbb{P}-$ complete probability space  $(\Omega,\mathcal{A},\mathbb{P})$. 
\end{corollary}
\nproof
From Theorem~\ref{thmpm}, this optimization problem is equivalent to the causal \textsc{Monge}-\textsc{Kantorovich} problem~(\ref{MKC1P}).  Hence, the result follows from \cite{CAUSAL}. \nqed

\begin{example}
For $m\in \mathbb{N}$, $\lambda >0$, denote by $\Gamma\left(m, \frac{1}{\lambda}\right)\in M_1(\mathbb{R})$ the \textsc{Borel} probability measure absolutely continuous with respect 
to the \textsc{Lebesgue} measure, whose density $\rho_{m,\frac{1}{\lambda}} :\mathbb{R} \to \mathbb{R}$ is given by $$\rho_{m,\frac{1}{\lambda}}(x) = \lambda \exp(-\lambda x) 
\frac{(\lambda x)^{m-1}}{(m-1)!} \mathds{1}_{[0,+\infty[}(x), \ \forall x\in \mathbb{R} ;$$ those are particular cases of the so-called \textit{gamma} distributions.
For  $n\in \mathbb{N}$, and $p \in \mathbb{N}$, define  $$\eta = \Gamma\left(n ,\frac{1}{\lambda}\right)$$ and $$\nu =\Gamma\left(n +p, \frac{1}{\lambda}\right),$$ for some $\lambda\in]0,+\infty[$. Further, let $(N_t)_{t\in\mathbb{R}_+}$ be a simple \textsc{Poisson} process with parameter $\lambda$, which is defined on some probability space $(\Omega,\mathcal{A},\mathbb{P})$ (for a definition, see \cite{CONT}, \cite{PROTTER} or \cite{STROOCK}), and set $$T_m= \inf(\{t \in \mathbb{R} : N_t \geq m\}),$$ which is a random variable on $(\Omega, \mathcal{A}, \mathbb{P})$, which describes the time of the $m-$th jump of the process, $\forall m\in \mathbb{N}$. Then, denoting by $p_{T_{n},T_{n+p}}$ the joint law of the jump times $T_n$ and $T_{n+p}$ on $(\Omega, \mathcal{A}, \mathbb{P})$, from Theorem~\ref{thmpm} we obtain $$p_{T_{n},T_{n+p}}\in \Pi_c(\eta,\nu),$$ so that  \textsc{Jensen}'s inequality ensures that it also attains the infimum of $$\inf_{\gamma \in \Pi_c(\eta,\nu)} \int_{\mathbb{R}\times \mathbb{R}} |x- y| \gamma(dx,dy).$$ 
\end{example}

\begin{remark} (On the identification of the optima, and the $c$-cyclical monotonicity.)
Let $\eta \in M_1(\mathbb{R})$, $\nu\in M_1(\mathbb{R})$, further assume that $c :\mathbb{R}\times \mathbb{R} \to \mathbb{R}$ is non-negative, continuous,  and that there exists $u\in \mathcal{L}^1(\eta)$, $v\in \mathcal{L}^1(\nu)$, such that $c(x,y) \leq u(x) +v(y)$, $\forall (x,y)\in \mathbb{R}\times \mathbb{R}$. Notice that,  from the proofs that the $c$- cyclical monotonicity is a necessary condition for optima of usual \textsc{Monge}-\textsc{Kantorovich} problem (see the fundamental theorem of optimal transport in \cite{AMB}),
 for $\gamma \in \Pi_c(\eta,\nu)$ to be an optimum of the causal mass transport problem~(\ref{MKC1P}), a necessary condition is that for any $n\in\mathbb{N}$, and for any $\{(x_i,y_i) : i\in\{1,...,n\}\}$ contained in the support of $\gamma$, such that \begin{equation} \label{esmple} \max_{i\in\{1,...,n\}} x_i < \min_{i\in\{1,...,n\}} y_i, \end{equation}  and for any permutation $\sigma : \{1,...,n\} \to \{1,...,n\}$, we necessarily have $$\sum_{i=1}^n c(x_i, y_i) \leq \sum_{i=1}^n c(x_i, y_{\sigma(i)}) ;$$ indeed, assuming~(\ref{esmple}), $(iii)$ of Lemma~\ref{empor} ensures that, following the proof by absurd which is proposed for the classical problem in \cite{AMB}, as soon as $\gamma \in \Pi_c(\eta,\nu)$, the plan $\widetilde{\gamma}$ of lower cost than $\gamma$ which is built in \cite{AMB} can be chosen to be an element of $\Pi_c(\eta,\nu)$, as far as, given $t_1\in ] \max_{i\in\{1,...,n\}} x_i ,\min_{i\in\{1,...,n\}} y_i[$, the neighbourhoods $(U_i)_{i\in\{1,...,n\}}$ (resp. $(V_i)_{i \in\{1,...,n\}}$) of the respective $(x_i)_{\in \{1,...,n\}}$ (resp. $(y_i)_{\in\{1,...,n\}}$),  where the mass may be modified, can be chosen to be such that $\bigcup_{i=1}^n U_i \subset ]-\infty, t_1[$ and  $\bigcup_{i=1}^n V_i .\subset]t_1,+\infty[$.
\end{remark}

\end{document}